\documentclass[12pt,leqno]{article}
\usepackage{amsmath, amsthm, amsfonts, amssymb, color}
\usepackage{mathrsfs}
\theoremstyle{definition}

\newcommand{\scr}[1]{\mathscr #1}
\definecolor{wco}{rgb}{0.5,0.2,0.3}

\numberwithin{equation}{section} \theoremstyle{remark}

\newcommand{\ua}{\uparrow}

\title{{\bf Analysis on Path Spaces over Riemmannian Manifolds with Boundary}\footnote{Supported in
 part by WIMCS, NNSFC(10721091) and the 973-Project.}
}
\author{
{\bf Feng-Yu Wang}\\
\footnotesize{School of Mathematical Sci. and Lab. Math. Com. Sys.,
Beijing Normal
University, Beijing 100875, China}\\
\footnotesize{and}\\ \footnotesize{Department of Mathematics,
Swansea University, Singleton Park, SA2 8PP, UK}\\ \footnotesize{Email: wangfy@bnu.edu.cn;
F.Y.Wang@swansea.ac.uk}}
\begin{document}
\def\R{\mathbb R}  \def\ff{\frac} \def\ss{\sqrt} \def\B{\mathbf
B}
\def\N{\mathbb N} \def\kk{\kappa} \def\m{{\bf m}}
\def\dd{\delta} \def\DD{\Delta} \def\vv{\varepsilon} \def\rr{\rho}
\def\<{\langle} \def\>{\rangle} \def\GG{\Gamma} \def\gg{\gamma}
  \def\nn{\nabla} \def\pp{\partial} \def\EE{\scr E}
\def\d{\text{\rm{d}}} \def\bb{\beta} \def\aa{\alpha} \def\D{\scr D}
  \def\si{\sigma} \def\ess{\text{\rm{ess}}}
\def\beg{\begin} \def\beq{\begin{equation}}  \def\F{\scr F}
\def\Ric{\text{\rm{Ric}}} \def\Hess{\text{\rm{Hess}}}
\def\e{\text{\rm{e}}} \def\ua{\underline a} \def\OO{\Omega}  \def\oo{\omega}
 \def\tt{\tilde} \def\Ric{\text{\rm{Ric}}}
\def\cut{\text{\rm{cut}}} \def\P{\mathbb P} \def\ifn{I_n(f^{\bigotimes n})}
\def\C{\scr C}      \def\aaa{\mathbf{r}}     \def\r{r}
\def\gap{\text{\rm{gap}}} \def\prr{\pi_{{\bf m},\varrho}}  \def\r{\mathbf r}
\def\Z{\mathbb Z} \def\vrr{\varrho} \def\ll{\lambda}
\def\L{\scr L}\def\Tt{\tt} \def\TT{\tt}\def\II{\mathbb I}
\def\i{{\rm i}}\def\Sect{{\rm Sect}}\def\E{\mathbb E} \def\H{\mathbb H}
\def\M{\scr M}\def\Q{\mathbb Q} \def\texto{\text{o}}

\maketitle
\begin{abstract} By using Hsu's multiplicative functional for the
Neumann heat equation, a natural damped gradient operator is defined
for the reflecting Brownian motion on compact manifolds with
boundary. This operator is  linked to   quasi-invariant flows in
terms of  a integration by parts formula, which leads to the
standard log-Sobolev inequality for the associated Dirichlet form on
the path space.
\end{abstract} \noindent

 \vskip 2cm

\section{Introduction}
Stochastic analysis on the path space over a complete Riemannian
manifold without boundary has been well developed since 1992 when B.
K. Driver \cite{D} proved the quasi-invariance theorem for the
Brownain motion on compact Riemannian manifolds. A key point of the
study is to first establish an integration by parts formula for the
associated gradient operator induced by the quasi-invariant flow,
then prove functional inequalities for the corresponding Dirichlet
form (see e.g. \cite{F,H,CHL} and references within). Moreover, some
efforts have been made for the study of geometry and topology on
Riemannian path or loop spaces (see e.g. \cite{EL} and references
within).

On the other hand, however,   the analysis on the path space over a
manifold with boundary is still very open. To see this, let us
mention  \cite{Z} where an integration by parts formula was
established on the path space of the one-dimensional reflecting
Brownian motion.  Let e.g. $X_t=|b_t|$, where $b_t$ is the
one-dimensional Brownian motion.  For   $h\in C([[0,T];\R)$ with
$h_0=0$ and  $\int_0^T |\dot h_t|^2\d t<\infty$,  let $\pp_h$ be the
derivative operator induced by the flow $X+\vv h$, i.e.

$$\pp_h F = \sum_{i=1}^n h_{t_i} \nn_i f(X_{t_1},\cdots, X_{t_n}),$$
where $n\in \mathbb N, 0<t_1<\cdots<t_n\le T$ and $F(X)=
f(X_{t_1},\cdots, X_{t_n})$ for some $f\in C^\infty(M^n).$  As the
main result of \cite{Z}, when $h\in C_0^2(0,T)$, \cite[Theorem
2.3]{Z} provides an integration by parts formula for $\pp_h$ by
using an infinite-dimensional generalized functional in the sense of
Schwartz.  Since for non-trivial $h$  the flow is not
quasi-invariant, this integration by parts formula can not be
formulated by using the distribution of $X$ with a density function,
and  the induced gradient operator does not provide a Dirichlet form
on the $L^2$-space of the distribution of $X$.

In this paper, we shall define  quasi-invariant flows on a
$d$-dimensional  Riemannian manifolds with boundary for all $h\in
\H$ in an intrinsic way, where

$$\H:= \bigg\{h\in C( [0,T]; \R^d):\ h_0=0, \int_0^T |\dot h_t|^2\d t<\infty\bigg\}$$
is the Cameron-Martin space. When $M$ is a half-space of $\R^d$,
which essentially reduces to the one-dimensional setting,
quasi-invariant flows has been constructed in \cite[\S 4(a)]{B} by
solving SDEs with reflecting boundary. We shall modify the idea  to
the reflecting Brownian motion on a manifold with boundary. By
establishing integration by parts formula, these flows will be
linked to a damped gradient operator defined by using Hsu's
multiplicative functionals constructed in \cite{H}. Form this we
will derive the Gross log-Sobolev inequality  for the associated
Dirichlet form.

To explain  the  idea of the study in a simple way, we first
consider the one-dimensional situation. Let $l_t$ be the local time
of $X_t:= |b_t|$ at point $0$. We have

$$\d X_t = \d b_t + \d l_t,\ \ X_0=0.$$ Now, for any $h\in \H$ and
$\vv>0$, let $X_t^{\vv,h}$ and its local time $l_t^{\vv,h}$ at $0$
solve the equation

$$ \d X_t^{\vv,h}= \d b_t +\vv \dot h_t\d t +\d l_t^{\vv,h},\ \
X_0^{\vv,h}=0.$$  By the Girsanov theorem $\{b_t+\vv h_t:\ 0\le t\le
T\}$ is a Brownian motion under the probability $R_\vv \P$, where

$$R_\vv =\exp\bigg[-\vv \int_0^T \dot h_t \d b_t -\ff {\vv^2} 2
\int_0^T |\dot h_t|^2 \d t\bigg]$$ is a functional of $X$ since $\d
b_t= \d X_t -\d l_t$. Thus, the distribution of $X^{\vv,h}$ under
$R_\vv \P$ coincides with that of $X$ under $\P$. Therefore,  the
flow $X^{\vv,h}$ is quasi-invariant. Moreover, it is easy to see
that
 $X_t^{\vv,h}= |b_t+\vv h_t|.$ So, for a cylindrical
 function $\gg\mapsto F(\gg)= f(\gg_{t_1},\cdots, \gg_{t_n})$, where $n\ge 1,
 0< t_1<\cdots <t_n\le T,$
 $f\in C_0^\infty(\R_+^n)$ and $\gg\in C([0,1]; [0,\infty))$ with
 $\gg_0=0,$ one has

 \beg{equation*}\beg{split} &\lim_{\vv\downarrow 0} \ff{F(X^{\vv,h})-
 F(X)}\vv =\sum_{i=1}^n \text{sgn} (b_{t_i}) h_{t_i} \pp_i
 f(X_{t_1},\cdots, X_{t_n})\\
 &=\sum_{i=1}^n \text{sgn}(X_{t_i}-
 L_{t_i}) h_{t_i} \pp_i f(X_{t_1},\cdots,
 X_{t_n}),\end{split}\end{equation*} which is a functional of $X$.
 Let $\tt f(x_1,\cdots, x_n)= f(|x_1|, \cdots, |x_n|)$. We have

\beq\label{1.0} D_h^0 F:=\lim_{\vv\downarrow 0} \ff{F(X^{\vv,h})-
 F(X)}\vv = \sum_{i=1}^n h_{t_i} \pp_i \tt f(b_{t_1}, \cdots,
 b_{t_n}).\end{equation} Combining this with the known integration by parts formula
 for the Brownian motion, we obtain

\beq\label{00}\E D_h^0 F= \E\bigg\{F(X) \int_0^T \dot h_t \d
b_t\bigg\}.\end{equation} Furthermore,
 let the gradient of $F$ be fixed as an $\H$-valued random variable such
 that $\<D^0F,h\>_\H= D_h^0 F, h\in \H.$ So,

 $$D^0F= \sum_{i=1}^n (t_i\land t) \text{sgn}(X_{t_i}-l_{t_i}) h_{t_i} \pp_i
 f(X_{t_1},\cdots, X_{t_n}).$$ Let $\mu$ be the distribution of
 $X$. By (\ref{00}), the    form

 $$\EE(F,G)=\E(D^0F, D^0G\>_\H$$ defined for cylindrical functions $F$ and $G$
 is closable in $L^2(\mu)$, and the closure $(\EE,\D(\EE))$ is a conservative
 Dirichlet form. Finally, by the known log-Sobolev inequality on the path space
 of the Brownain motion, this Dirichlet form  satisfies the log-Sobolev
 inequality

 $$\mu(F^2\log F^2)\le 2\EE(F,F),\ \ F\in \D(\EE).$$

 The main purpose of this paper is to realized the above idea on the path space of the reflecting
 Brownian motion on a Riemannian manifold with boundary. In this
 case we no longer have explicit expression of $D^0$. But in Section
 2 we shall present an integration by parts formula, which identifies
 the adapted projection of $D^0$ and that of the damped gradient
 operator induced by Hsu's multiplicative functional constructed in
 \cite{H}. this
 integration by parts formula will be proved in Section 3. Finally, using the resulting integration by
 parts formula, the standard log-Sobolev
 inequality will be addressed in Section 4.

\section{Damped Gradient and Integration by Parts}

Let $M$ be a $d$-dimensional compact connected Riemannian manifold
with boundary $\pp M$. Let $o\in M$ and $T>0$ be  fixed. Then the
path space for the reflecting Brownian motion on $M$ starting at $o$
is

$$W= \{\gg\in C([0,T]; M):\ \gg_0=o\}.$$ Let $B_t$ be the
$d$-dimensional Brownian motion on a complete probability space
$(\OO,\F,\P)$ with natural  filtration $\{\F_t\}_{t\ge 0}.$ For any
$x\in M$, let $O_{x}M$ be the set of all orthonormal bases for the
tangent space $T_{x}M$ at point $x$, and let $O(M):= \cup_{x\in M}
O_x(M)$ be the frame bundle. Then for any $X_0\in M$, the reflecting
Brownian motion can be constructed by solving the SDE

\beq\label{2.1} \d X_t= u_t\circ \d B_t +N(X_t)\d l_t,\end{equation}
where $u_t\in O_{X_t}(M)$  is the horizontal lift of $X_t$ on the
frame bundle $O(M)$, and $l_t$ is the local time of $X_t$ on the
boundary $\pp M.$ Let $\mu$ be the distribution of $X:=\{X_t:\ 0\le
t\le T\}$ for $X_0=o$. Then $\mu$  is   a probability measure on the
path space $W$.

To define the damped gradient operator, let us introduce the
multiplicative functional constructed in  \cite{H}. To this end, we
need to introduce some $\R^d\bigotimes\R^d$-valued functionals on
the frame bundle. Let $\Ric$ be the Ricci curvature on $M$ and $\II$
the second fundamental form on $\pp M$. For any $u\in O(M)$, let

$$R_u (a,b)= \Ric(au, bu),\ \  a,b\in \R^d.$$ Let $\pi_\pp: TM\to T\pp M$ be the orthogonal projection
at points on $\pp M$, and let   $\pi: O(M)\to M$ is the canonical
projection. For any $u\in O(M)$ with $\pi u\in\pp M$, let

$$\II_u(a,b)= \II(\pi_\pp au, \pi_\pp bu),\ \ a,b \in \R^d.$$ Finally,
let $N$ be the inward unit normal vector field on $\pp M$. For $\pi
u\in \pp M$, let

$$P_u(a,b)= \<ua, N\>\<ub, N\>,\ \ a,b\in \R^d.$$

For any $u_0\in O(M)$, let $X_t$ be the reflecting Brownian motion
on $M$ with horizontal lift $u_t$. For any $\vv>0$, let $  Q_t^\vv$
solve the following SDE on $\R^d\bigotimes\R^d$:

\beq\label{PP1} \d Q_t^\vv =- Q_t^\vv \Big\{\ff 1 2 R_{u_t} \d t +
\big(\vv^{-1}P_{u_t} +\II_{u_t}\big)\d l_t\Big\},\ \ \ Q_0^\vv =
I.\end{equation} According to \cite[Theorem 3.4]{H}, when
$\vv\downarrow 0$ the process $Q_t^\vv$ converges in $L^2$ to an
adapted right-continuous process $Q_t$ with left limit, such that
$Q_tP_{u_t}=0$ if $X_t=\pi u_t\in \pp M.$ Consequently, if $\Ric\ge
-K$ and $\II\ge -\si$ for some continuous functions $K$ and $\si$ on
$M$, then

$$\|Q_t\|\le \exp\bigg\{\ff 1 2 \int_0^t K(X_s)\d s + \int_0^t
\si(X_s)\d l_s\bigg\},\ \ t\ge 0,$$ where $\|\cdot\|$ is the
operator norm on $\R^d$. In particular, $\E\|Q_t\|^p<\infty$ holds
for any $p>1$. For $f\in C(M)$, let

$$P_t f(x)= \E f(X_t^x),\ \ \ x\in M,$$ where and in the sequel, $X_t^x$ denotes the solution to (\ref{2.1}) for
$X_0=x$. Then $P_t$ is  the Neumann semigroup.
  By \cite[Theorem 4.2]{H}
(see also the last display in the proof of \cite[Theorem 5.1]{H}),
 $s\mapsto Q_s u_s^{-1} \nn P_{t-s} f (X_s)$ is a martingale. So,

\beq\label{2.2} u_0^{-1}  \nn P_t f(x)  = \E \big\{  Q_t^x
(u_t^x)^{-1} \nn f(X_t^x)\big\},\ \ x\in M, u_0\in
O_x(M),\end{equation} where  $Q^x_t$ and  $u^x_t$ are the
multiplicative functional and horizontal lift of $X_t^x$.

In general, for $s\ge 0$, let $( Q_{s, t+s})_{t\ge 0}$ be the
associated multiplicative functional for the process
$(X_{t+s})_{t\ge 0}.$ We have

\beq\label{Q} Q_{s,t} Q_{t,r}= Q_{s,r},\ \ 0\le s\le t\le
r.\end{equation} We shall use these multiplicative functionals to
define the damped gradient operator (see \cite{FM} for the damped
gradient operator for manifolds without boundary).

Let

$$ \scr FC^\infty=\big\{ W\ni \gg\mapsto f(\gg_{t_1},\cdots, \gg_{t_n}): n\ge 1,
 0<t_1<\cdots< t_n\le T, f\in C^\infty(M)\big\}$$ be the class of
 smooth cylindrical functions on $W$. For any $F\in
\scr FC^\infty$ with $F(\gg)=f(\gg_{t_1},\cdots, \gg_{t_n}),$ define
the damped gradient $DF$ as an $\H$-valued random variable by
setting $(DF)_0=0$ and

$$\ff{\d}{\d t}(DF)_t= \sum_{i=1}^n 1_{\{t<t_i\}}  Q_{t,   t_i} u_{t_i}^{-1}
\nn_i f(X_{t_1},\cdots, X_{t_n}),\ \ t\in [0,T],$$ where $\nn_i$
denotes the gradient operator w.r.t. the $i$-th component. Then, for
any $\H$-valued random variable $h$, we have

\beq\label{Dh} D_h F:= \<DF, h\>_\H= \sum_{i=1}^n \int_0^{t_i}
\<u_{t_i}^{-1} \nn_i f(X_{t_1},\cdots, X_{t_n}), Q_{t, t_i}^* \dot
h_t\>\d t.\end{equation}

Next, let $\tt \H$ denote the set of all square-integrable
$\H$-valued adapted random variables, i.e.

$$\tt \H=\big\{h\in L^2(\OO\to \H;\P):\ h_t \ \text{is}\ \F_t\text{-measurable},\ t\in [0,T]\big\}.$$
 Then $\tt\H$ is a Hilbert
space with  inner product

$$\<h,\tt h\>_{\tt \H}:= \E \int_0^T \<\dot h_t, \dot{\tt h}_t\>\d
t=\E\<h,\tt h\>_\H,\ \ h,\tt h\in \tt \H.$$

To describe $DF$ by using a quasi-invariant flow, for $h\in\tt\H$
and $\vv>0$ let $X_t^{\vv,h}$ solve the SDE

\beq\label{2.3} \d X_t^{\vv,h} = u_t^{\vv,h} \circ \d B_t
+N(X_t^{\vv,h})\d l_t^{\vv,h} +\vv \dot h_t u_t^{\vv,h}\d t,\ \
X_0^{\vv,h}=X_0=\pi u_0,\end{equation} where $l_t^{\vv,h}$ and
$u_t^{\vv,h}$ are, respectively, the local time on $\pp M$ and the
horizontal lift on $O(M)$ for $X_t^{\vv,h}$. To see that
$\{X^{\vv,h}\}_{\vv\ge 0}$ has the flow property, let

$$\Theta: W_0:=\{\oo\in C([0,T]):\ \oo_0=0\}\to W$$ be measurable such
that $X=\Theta(B)$. For any $\vv>0$ and a function $\Phi: W_0\to W$,
let $(\theta_\vv^h \Phi)(\oo)= \Phi(\oo+\vv h).$ Then $X^{\vv,h}=
(\theta_\vv^h \Theta)(B), \vv\ge 0$. Hence,

$$X^{\vv_1+\vv_2,h}= \theta_{\vv_1}^h X^{\vv_2,h},\ \ \
\vv_1,\vv_2\ge 0.$$ We shall try to link the multiplicative
functional $ Q$ to the vector field (if exists) generating the flow
$X^{\vv,h}.$

First of all, let us explain that the flow $X^{\vv,h}$ is
quasi-invariant. Let

$$R^{\vv,h}= \exp\bigg[\vv \int_0^T\<\dot h_t, \d B_t\> -\ff {\vv^2}
2 \int_0^T |\dot h_t|^2\d t\bigg].$$ By the Girsanov theorem,

$$B_t^{\vv,h}:= B_t-\vv h_t$$ is the $d$-dimensional Brownian motion
under the probability $R^{\vv,h}\P.$ Thus, the distribution of $X$
under $R^{\vv,h}\P$ coincides with that of $X^{\vv,h}$ under $\P.$
Therefore,   $X^{\vv,h}$ is quasi-invariant.

The following   integration by parts formula provides a link between
 the damped gradient $D$ and the flow $X^{\vv,h}$.

\beg{thm}\label{T2.1} For any $u_0\in O(M)$ and $F\in \scr
FC^\infty$,

$$ \E\big\{D_h F\big\}= \lim_{\vv\downarrow 0} \E
\ff{F(X^{\vv,h})-F(X)}\vv =\E \bigg\{F(X)\int_0^T \<\dot h_t, \d
B_t\>\bigg\}$$ holds for all $h\in \tt \H_b$, the set of all
elements in $\tt \H$ with bounded $\|h\|_\H.$ \end{thm}

\paragraph{Remark 2.1.} Since $\tt \H_b$ is dense in $\tt \H$, the above result  implies
that the projection of $D$ onto $\tt\H$ can be determined  by  the
flows $X^{\vv,h}, h\in \tt\H_b.$ But in general

\beq\label{Q*}D_h F= \lim_{\vv\downarrow 0} \ff{F(X^{\vv,h})-F(X)}
\vv,\ \ h\in\tt \H\end{equation} does not hold, so that    the flow
$\{X^{\vv,h}\}$ is not generated by the vector field

$$W\ni \gg\mapsto \bigg\{u_t(\gg) \int_0^t  Q_{s,t}^*(\gg) \dot h_s \d s\in T_{\gg_t}M: 0\le t\le T\bigg\},$$
where $u_s(\gg)$ and $ Q_s(\gg)$ are the horizontal lift and the
multiplicative functional of $\gg$ respectively.

To disprove (\ref{Q*}), let us consider $M=[0,1]\subset \R$ and
$X_0=0$. Let $h_t=t$ and $F(\gg)= gg_1.$ By (\ref{1.0}) we have

\beq\label{P1} \lim_{\vv\downarrow 0} \ff{F(X^{\vv,h})- F(X)}\vv
=f'(X_1)\,\text{sgn}(X_1-l_1)= \text{sgn}(X_1-l_1),\end{equation}
provided

$$\tau_1:= \inf\{t>0:\ X_t=1\} >1.$$ On the other hand, for the one-dimensional case we
have $R_u=\II_u=0$ and $P_u=1.$ Then by (\ref{PP1})

$$Q_{s,t}^\vv= \exp\big[-\vv^{-1} (l_t-l_s)\big]>0,\ \ s\le t.$$
This implies that  $Q_{t,1}\ge 0$ for all $t\in [0,1]$. Combining
this with (\ref{2.3}) we obtain

\beq\label{P2} D_h F=   f'(X_1)  \int_0^1 Q_{t,1}\d t=\int_0^1
Q_{t,1}\d t\ge 0.\end{equation}  Since

$$\P\big( X_1-l_1<0,\tau_1>1\big)= \P\Big(\inf_{s\in [0,1]} B_s<0,
\sup_{s\in [0,1]}B_s<1\Big)>0,$$ where $B_s$ is now the
one-dimensional Brownian motion. Combining this with    (\ref{P1})
and (\ref{P2}), we see that   $(\ref{Q*})$ does not hold.

\

To prove Theorem \ref{T2.1}, we need some preparations. In
particular, we shall use (\ref{2.2}) and a conducting argument as in
\cite{H2} for the case without boundary.

\section{Proof of Theorem \ref{T2.1}}

\beg{lem}\label{L3.1} Let $u_0\in O(M)$ and $F\in \F C^\infty.$ Then

$$\lim_{\vv\downarrow 0} \E \ff{F(X^{\vv,h})-F(X)}\vv =\E
\bigg\{F(X)\int_0^T \<\dot h_t, \d B_t\>\bigg\}$$ holds for all
$h\in \tt \H_b$.\end{lem}

\beg{proof} Let $B_t^{\vv,h}= B_t-\vv h_t,$ which is the
$d$-dimensional Brownian motion under $R^{\vv,h}\P.$ Reformulating
 (\ref{2.1}) as

 $$\d X_t= u_t\circ \d B_t^{\vv,h} + N(X_t)\d l_t +\vv \dot h_t
 u_t\d t,\ \ X_0=\pi u_0.$$ By the weak uniqueness of (\ref{2.3}),
 we conclude that the distribution of $X$ under $R^{\vv,h}\P$
 coincides with that of $X^{\vv,h}$ under $\P$. In particular, $\E
 F(X^{\vv,h})= \E [R^{\vv,h} F(X)]$. Thus,

\beg{equation*}\beg{split} &\lim_{\vv\downarrow 0} \E
\ff{F(X^{\vv,h})-F(X)}\vv
 =\lim_{\vv\downarrow 0} \E\Big\{F(X)
 \cdot\ff{R^{\vv,h}-1}\vv\Big\}\\
 &= \E\bigg\{F(X) \int_0^T \<\dot h_t, \d
 B_t\>\bigg\},\end{split}\end{equation*} where the last step is due
 to the dominated convergence theorem since $\{R^{\vv,h}\}_{\vv\in
 [0,1]}$ is uniformly integrable for $h\in \tt\H_b$. \end{proof}

\beg{lem}\label{L3.2} For any $n\ge 1, 0<t_1<\cdots<t_n\le T$, and
 $f\in C^\infty(M^n)$,

$$u_0^{-1}  \nn_x \E f(X_{t_1}^x,\cdots, X_{t_n}^x)
 =  \sum_{i=1}^n \E \big\{M^x_{t_i} (u_{t_i}^x)^{-1} \nn_i
f(X_{t_1}^x,\cdots, X_{t_n}^x)\big\} $$ holds for all $x\in M$ and
$u_0\in O_x(M),$ where $\nn_x$ denotes the gradient w.r.t.
$x$.\end{lem}

\beg{proof} By (\ref{2.2}), the desired assertion holds for $n=1$.
Assume that it holds for $n=k$ for some natural number $k\ge 1$. It
remains to prove the assertion for $n=k+1.$ To this end, set

$$g(x)= \E  f(x, X_{t_2-t_1}^x,\cdots, X_{t_{k+1}-t_1}^x),\ \ x\in M.$$
By the assumption for $n=k$ we have

$$u_0^{-1} \nn g(x) =    \sum_{i=1}^{k+1}
\E \big\{M^x_{t_i-t_1}(u_{t_i-t_1}^x)^{-1} \nn_i f(x,
X_{t_2-t_1}^x,\cdots, X_{t_{k+1}-t_1}^x)\big\}$$ for all $ x\in M,
u_0\in O_x(M).$   Combining this with the assertion for $k=1$ and
using the Markov property, we obtain

\beg{equation*}\beg{split}& u_0^{-1} \nn_x  \E  f(X_{t_1}^x,\cdots,
X_{t_{k+1}}^x)   =u_0^{-1} \nn_x
\E g(X_{t_1}^x)  \\
&= \E   \big\{Q_{t_1}^x (u_{t_1}^x)^{-1} \nn g(X_{t_1}^x)\big\}  =
\sum_{i=1}^{k+1} \E \big\{ Q_{t_i}^x (u_{t_i}^x)^{-1} \nn_i
f(X_{t_1}^x,\cdots, X_{t_{k+1}}^x)\big\}.\end{split}\end{equation*}
\end{proof}

\beg{lem}\label{L3.3} Let $f\in C^\infty(M)$. Then for any
 $u_0\in O(M)$ and $t>0$,

$$\E\bigg\{f(X_{t})\int_0^{t} \<\dot h_s, \d B_s\>\bigg\}= \E \int_0^{t}\<
 u_{t}^{-1} \nn f(X_{t}),  Q_{s,t}^*\dot h_{s}\>\d s,\ \
h\in\tt\H, t_1\in [0,T].$$\end{lem}

\beg{proof} Noting that

$$\ff{\d}{\d s} P_s f= \ff 1 2\DD P_s f,\ \ NP_sf|_{\pp M}=0,\ \
s>0,$$ by (\ref{2.1}) and the
 It\^o formula we obtain

$$\d (P_{t-s}f)(X_s)= \<\nn P_{t-s} f(X_s), u_s\d B_s\>,\ \ s\in
[0,t).$$ This implies

$$f(X_{t})= P_t f(X_0) +\int_0^{t} \<u_s^{-1}\nn P_{t-s}f(X_s), \d
B_s\>,\ \ s\in [0,t].$$ Therefore,

\beg{equation}\label{BB} \E\bigg\{f(X_{t}) \int_0^{t} \<\dot h_s,\d
B_s\>\bigg\}= \E\int_0^{t}\<u_s^{-1} \nn P_{t-s} f(X_s), \dot
h_s\>\d s.\end{equation}   By (\ref{2.2}) and the Markov property we
have

$$u_s^{-1} \nn P_{t-s} f(X_s)
=\E\big( Q_{s,t}u_{t}^{-1}\nn f(X_{t})\big|\F_{s}\big).$$ So, the
desired formula follows from (\ref{BB}) since $\dot h_s$ is
$\F_s$-measurable.
\end{proof}

 As a   consequence of (\ref{2.2}) and Lemma \ref{L3.3},  we have the
following Bismut formula.

\beg{cor} Let $P_t f(x)= \E^x f(X_t),\ t\ge 0, x\in M, f\in C(M).$
Then for any $v\in T_x M$ and any $h\in \tt H$ with $h_t= u_0^{-1}
v$,

$$\<v, \nn P_t f(x)\>=  \E^x \bigg\{f(X_t) \int_0^t
\< Q_s^*\dot h_s,\d B_s\>\bigg\}.$$\end{cor}

\beg{proof} By (\ref{Q}) and applying Lemma \ref{L3.3} to $\tt
h\in\tt \H$ in place of $h$, where $\dot{\tt h}_s= Q_s^* \dot h_s,$
we obtain

$$\E\bigg\{f(X_t) \int_0^t \<Q_s^* \dot h_s, \d B_s\>\bigg\} =\E \int_0^t
\<u_t^{-1} \nn f(X_t), Q_t^* \dot h_s\>\d s= \E \<Q_t u_t^{-1}\nn
f(X_t), u_0^{-1}v\>.$$ Then the proof is completed by combining this
with (\ref{2.2}). \end{proof}

\ \newline \emph{Proof of Theorem \ref{T2.1}.} By Lemma \ref{L3.1},
it suffices to prove

\beq\label{3.*} \E\{D_h F\} =\E \bigg\{F(X)\int_0^T\<\dot h_t, \d
B_t\>\bigg\},\ \ h\in\tt \H\end{equation} for $F = f(X_{t_1},
\cdots, X_{t_n})$ with $f\in C^\infty(M^n),$ where $n\ge 1,
0<t_1<\cdots<t_n\le T$.  According to Lemma \ref{L3.3}, (\ref{3.*})
holds for $n=1.$ Assuming (\ref{3.*}) holds for $n=k$ for some $k\ge
1$, we aim to prove it for $n=k+1.$ To this end, let

$$g(x)= \E  f(x, X_{t_2-t_1}^x,\cdots, X_{t_{k+1}-t_1}^x),\ \ x\in M.$$
By the result for $n=1$ and the Markov property,

\beq\label{2.4} \beg{split}  \int_0^{t_1} \E\< u_{t_1}^{-1}\nn
g(X_{t_1}),  Q_{t,t_1}^*\dot h_{t}\>\d t  &= \E\bigg\{\E(F(X)|\scr
F_{t_1})\int_0^{t_1}\<\dot h_t,\d
B_t\>\bigg\}\\
&= \E\bigg\{F(X)\int_0^{t_1} \<\dot h_t, \d
B_t\>\bigg\}.\end{split}\end{equation} On the other hand, by
(\ref{Q}), Lemma \ref{L3.2} and  the Markov property,

\beg{equation*}\beg{split} &\int_0^{t_1}\E\< u_{t_1}^{-1}\nn
g(X_{t_1}),  Q_{t,t_1}^* \dot h_{t}\>\d t\\ &=
\int_0^{t_1}\E\Big\<\E\Big(\sum_{i=1}^{k+1}  Q_{t_1,t_i}
u_{t_i}^{-1} \nn_i f(X_{t_1},\cdots, X_{t_{k+1}})\Big|\scr
F_{t_1}\Big),Q_{t,t_1}^*\dot h_t \Big\>\d t\\
&= \sum_{i=1}^{k+1}\int_0^{t_1}\< u_{t_i}^{-1} \nn_i
f(X_{t_1},\cdots, X_{t_{k+1}}),
  Q_{t, t_i}^*\dot h_t \>\d t.\end{split}\end{equation*} Combining this with
(\ref{Dh}) and (\ref{2.4}) we obtain

\beq\label{2.5} \beg{split} \E\big\{D_h F\big\} =
&\E\bigg\{F(X)\int_0^{t_1} \<\dot h_t, \d B_t\> \bigg\}\\& +\E
\sum_{i=2}^{k+1} \int_{t_1}^{t_i} \< u_{t_i}^{-1} \nn_i f(X_{t_1},
\cdots, X_{t_{k+1}}),  Q_{t, t_i}^*\dot h_t \>\d t
.\end{split}\end{equation}  By the Markov property and the
assumption for $n=k$, we have

$$\sum_{i=2}^{k+1} \E\int_{t_1}^{t_i} \< u_{t_i}^{-1} \nn_i
f(X_{t_1}, \cdots, X_{t_{k+1}}),   Q_{t, t_i}^*\dot h_t \>\d t  =\E
\bigg\{F(X)\int_{t_1}^{T}\<\dot h_t,\d B_t\>\bigg\}.$$Combining this
with (\ref{2.5}) we complete the proof.\qed

\section{The Log-Sobolev Inequality}

Let $\mu$ be the distribution of $X$ with $X_0=o$, and let

$$\EE (F,G)= \E\<DF, DG\>_\H,\ \ F,G\in \scr FC^\infty.$$
 Since both $DF$ and $DG$ are functionals of $X$,
$(\EE,\scr FC^\infty)$ is a positive  bilinear form on $L^2(W;\mu)$.
It is standard that the integration by parts formula (\ref{3.*})
implies the closability of the form (see Lemma \ref{L4.0}).  We
shall use $(\EE,\D(\EE))$ to denote the closure of $(\EE,\scr
FC^\infty)$. Moreover, (\ref{3.*})  also implies    the Clark-Ocone
type martingale representation formula (see Lemma \ref{L4.1}), which
leads to the standard Gross \cite{G} log-Sobolev inequality. It is
well known that the log-Sobolev inequality   implies that the
associated Markov semigroup is hypercontractive and converges
exponentially to $\mu$ in the sense of relative entropy.

\beg{lem}\label{L4.0} $(\EE,\scr FC^\infty)$ is closable in
$L^2(W;\mu)$.
\end{lem}

\beg{proof} Although the proof is standard by using the integration
by parts formula, we include it here for completeness.  Let
$\{F_n\}_{n\ge 1}\subset \scr FC^\infty$ such that $\EE(F_n,F_n)\le
1$ for all $n\ge 0$ and $\mu(F_n^2)+\EE(F_n-F_m,F_n-F_m)\to 0$ as
$n,m\to\infty.$ We aim to prove that $\EE(F_n,F_n)\to 0$ as
$n\to\infty.$ Since

$$\EE(F_n,F_n)= \EE(F_n,F_n-F_m)+ \EE(F_n,F_m)\le \ss{\EE(F_n-F_m,
F_n-F_m)} +\EE(F_n,F_m),$$ it suffices to show that for any $G\in
\scr FC^\infty$, one has $\EE(F_n,G)\to 0$ as $n\to\infty$. To this
end, let $\{h^i\}_{i\ge 1}$ be an ONB on $\H$. For any $\vv>0$ there
exists $k\ge 1$ such that

$$\Big|\EE(F_n,G)-\sum_{i=1}^k  \E(D_{h^i} F_n)(D_{h^i}G) \Big|<\vv,$$ where
$D_hF:= \<DF,h\>_{\H}$ for $F\in \scr FC^\infty$ and $h\in\H.$ Since
$\scr FC^\infty$ is dense in $L^2(W;\mu)$, there exists $G_i\in \scr
FC^\infty$ such that

$$\E|D_{h^i} G- G_i|^2<\vv,\ \ 1\le i\le k.$$ Therefore,

$$|\EE(F_n,G)|\le 2\vv + \sum_{i=1}^k \big|\E \<G_iDF_n, h_i\>_\H\big|.$$
Noting that $G_i DF_n= D(F_nG_i)- F_nDG_i$, by (\ref{3.*}) we obtain

$$|\EE(F_n,G)| \le 2\vv +\sum_{i=1}^k \bigg|\E
\bigg[F_n(X)\bigg\{G_i(X)\int_0^T \<\dot{h}^i_t, \d B_t\> -D_{h^i}
G_i\bigg\}\bigg]\bigg|.$$ Since $\mu(F_n^2)\to 0$ as $n\to\infty$,
by letting first $n\to\infty$ then $\vv\to 0$ we complete the proof.
\end{proof}

\beg{lem}\label{L4.1} For any $F\in \scr FC^\infty$, let $\tt D F$
be the projection of $DF$ on $\tt\H$, i.e.

$$\ff{\d}{\d t}  (\tt DF)_t = \E\Big(\ff{\d }{\d t} (DF)_t\Big|\F_t\Big),\ \ t\in [0,T], (\tt DF)_0=0.$$
Then

$$F(X)= \E F(X) +\int_0^T \Big\<\ff{\d }{\d t} (\tt DF)_t, \d
B_t\Big\>.$$\end{lem}

\beg{proof}   By Theorem \ref{T2.1}, we have

\beq\label{4.1} \E\<h, \tt DF\>_\H= \E\bigg\{F(X)\int_0^T \<\dot
h_t, \d B_t\>\bigg\},\ \ h\in\tt\H.\end{equation}  On the other
hand, by the martingale representation, there exists a predictable
process $\bb_t$ such that

\beq\label{4.2} \E(F(X)|\F_t) = \E F(X) + \int_0^t \<\bb_s, \d
B_s\>,\ \ \ t\in [0,T].\end{equation} Let

$$\varphi_t= \int_0^t \bb_s\d s,\ \ \ t\in [0,T].$$ We have $\varphi\in
\tt\H$ and by (\ref{4.2}),

 $$\E \<h,\varphi\>_\H= \E\int_0^T \<\dot h_t,\bb_t\>\d t=
 \E\bigg\{F(X)\int_0^T \<\dot h_t,\d B_t\>\bigg\}$$ holds for all $h\in
 \tt\H$. Combining this with (\ref{4.1}) we conclude that $\tt
 DF=\varphi$. Therefore, the desired formula follows from (\ref{4.2}).
 \end{proof}

 \beg{thm} For any $T>0$ and   $o\in M$, $(\EE,\D(\EE))$ satisfies
 the following log-Sobolev inequality

 $$\mu(F^2\log F^2) \le 2 \EE(F,F),\ \ F\in \D(\EE),\
 \mu(F^2)=1.$$\end{thm}

 \beg{proof}  It suffices to prove for $F\in \F C^\infty$. Let $m_t=\E(F(X)^2|\F_t),\ t\in [0,T].$
By Lemma \ref{L4.1} and the It\^o formula,

$$\d m_t\log m_t= (1+\log m_t)\d m_t +\ff{|\ff{\d}{\d t} (\tt DF^2)_t|^2}{2m_t}\d t.$$ Thus,

\beg{equation*}\beg{split} \mu(F^2\log F^2) &= \E m_T\log m_T
=\int_0^T \ff{2 \E(F(X)
\ff{\d }{\d t} (DF)_t|\F_t)^2}{\E(F(X)^2|\F_t)}\d t\\
&\le 2 \int_0^T \E\Big| \ff{\d}{\d t} (DF)_t\Big|^2\d t= 2
\E\|DF\|_\H^2= 2\EE(F,F).\end{split}\end{equation*}\end{proof}

 \beg{thebibliography}{99}

\bibitem{B} J.-M. Bismut, \emph{the calculus of boundary processes,}
Ann. Sci. \'Ecole Norm. Sup. 17(1984), 507--622.

\bibitem{CHL} M. Capitaine, E. P. Hsu and M. Ledoux, \emph{Martingale
representation and a simple proof of logarithmic Sobolev
inequalities on path spaces,} Electron Comm. Prob. 2(1997) 71-81.

\bibitem{D} B. K. Driver, \emph{A Cameron-Martin type
quasi-invariance theorem for Brownian motion on a compact Riemannian
manifold,} J. Funct. Anal. 110(1992), 272--376.

\bibitem{EL} K. D.  Elworthy, K. D.,  X.-M. Li, \emph{An $L^2$ theory for differential forms on path
spaces I,}  J. Funct. Anal. 254 (2008),   196--245.

\bibitem{F} S. Fang, \emph{Inegalit\'e du type de Poinca\'e sur
l'espace des chemins riemanninens,} C.R.Acad.Sci.Paris 318(1994),
257--260.

\bibitem{G} L. Gross,  \emph{Logarithmic Sobolev inequalities,} Amer. J. Math.
 97(1976), 1061--1083.

  \bibitem{FM} S. Fang and P. Malliavin, \emph{Stochastic analysis on the path
space of a Riemannian manifold},  J. Funct. Anal. 118 (1993),
p.249--274.

\bibitem{H2} E. P. Hsu, \emph{Logarithmic-Sobolev inequalities on
path spaces over Riemannian manifolds,} Comm. Math. Phys. 189(1997),
9--16.

\bibitem{H} E. P. Hsu, \emph{Multiplicative functional for the heat
equation on manifolds with boundary,} Michigan Math. J. 50(2002),
351--367.

\bibitem{Z} L. Zambotti, \emph{Integration by parts on the law of the reflecting
Brownian motion,} J. Funct. Anal. 223(2005), 147-178.
 \end{thebibliography}

\end{document}